\documentclass[12pt]{article}

\usepackage{amssymb,amsfonts,amsmath}
\usepackage{tikz}

\makeatletter\@addtoreset{figure}{section}\makeatother

\makeatletter
\long\def\@makecaption#1#2{%
   \vskip 10\p@
   \setbox\@tempboxa\hbox{{#1} \ \ #2}%
   \ifdim \wd\@tempboxa >\hsize
       {#1}\ \ #2\par
   \else
       \hbox to\hsize{\hfil\box\@tempboxa\hfil}%
   \fi}
\makeatother

\newtheorem{theo}{Theorem}

\newtheorem{lem}[theo]{Lemma}
\newtheorem{exam}{Example}[section]

\makeatletter \@addtoreset{equation}{section}
\@addtoreset{theo}{section} \makeatother

\def\qed{\hfill \rule{4pt}{7pt}}
\def\pf{\noindent {\it Proof.} }

\title{A combinatorial proof of the Alladi-Gordon key identity for Schur's partition theorem}
\author{James J.Y. Zhao\\
\small Dongling School of Economics and Management,\\
\small University of Science and Technology Beijing, Beijing 100083, P.R. China\\
\small \texttt{zhaojy@ustb.edu.cn}}

\date{\small{Feb 6, 2012}}

\begin{document}

\maketitle

\begin{abstract}
The Alladi-Gordon identity plays an important role for the Alladi-Gordon generalization of Schur's partition theorem. By using Joichi-Stanton's insertion algorithm, we present an overpartition interpretation for the Alladi-Gordon key identity. Based on this interpretation, we further obtain a combinatorial proof of the Alladi-Gordon key identity by establishing an involution on the underlying set of overpartitions.
\end{abstract}


\section{Introduction}

Let $\mathbb{N}$ be the set of nonnegative integers. Let
\begin{eqnarray*}
(a)_k=(a;q)_k=\left\{
\begin{array}{@{}ll@{}l}
 (1-a)(1-aq)\cdots(1-aq^{k-1}),&\quad &{\rm if}\ k>0,\\
 \ 1,&\quad & {\rm if}\ k=0,
\end{array}\right.
\end{eqnarray*}
denote the common notation of $q$-shifted factorials \cite{GaRa90}. Given $j,k\in \mathbb{N}$, let
\begin{equation*}
{j\brack k}=\frac{(q;q)_j}{(q;q)_k(q;q)_{j-k}},
\end{equation*}
denote the Gaussian coefficients, which are also called as the $q$-binomial coefficients, or the
Gaussian polynomials \cite{Andrews76}. The main objective of this paper is to give a combinatorial proof of the following identity:
\begin{equation}\label{eq-ij2}
\sum\limits_{k=0}^{j}(q^{i-k+1};q)_k\,{j\brack k} q^{(i-k)(j-k)}=1,
\end{equation}
which we call \emph{the Alladi-Gordon key identity}, since it was first introduced by Alladi and Gordon \cite{Alladi-Gordon-93} in an equivalent form for the study of some generalization of Schur's celebrated partition theorem
of 1926.

Schur \cite{Schur26} proved that
the number of partitions of
$m$ into parts with minimal difference $3$ and with no consecutive
multiples of $3$  is equal to  the number of partitions of $m$ into distinct parts
$\equiv1,2\,(\!\!\!\mod 3)$. This significant result is now known as {\it  Schur's celebrated partition theorem
of 1926}. There are many proofs, refinements, and generalizations of
Schur's partition theorem, see \cite{Gleissberg, Gollnitz-67, Andrews67,
Andrews68, Andrews68-2, Andrews69, Andrews71, Bressoud80,
Bessenrodt91, Alladi-Gordon-93, Andrews94, Alladi-Gordon-95, AB2002,
Boulet-Pak, CorLov2006, Pak2006}.

From the viewpoint of generating functions, each partition theorem implies a corresponding $q$-identity.
The Alladi-Gordon key identity \eqref{eq-ij2} is essentially equivalent to the following $q$-identity \cite[Lemma 2]{Alladi-Gordon-93} corresponding to Alladi and Gordon's notable generalization of Schur's partition theorem,
\begin{equation}\label{AG93-Lem2}
\sum\limits_{0\leq m \leq min\{i,j\}} \frac{q^{T_{i+j-m}+T_m}}
{(q)_{i-m}(q)_{j-m}(q)_m}=\frac{q^{T_i+T_j}}{(q)_i(q)_j},
\end{equation}
where $i$ and $j$ are given nonnegative integers and $T_i=i(i+1)/2$ is
the $i$-th triangle number.

The Alladi-Gordon key identity turned out to have many interesting applications in
the theory of partitions. Alladi and Berkovich \cite[Eq.\,\,(2.1)]{AB2002}) obtained a double bounded version of Schur's
partition theorem by generalizing an equivalent form of \eqref{AG93-Lem2}. Alladi, Andrews and Gordon
\cite[Lemma 2]{AAG-95} introduced a three parameter generalization of
\eqref{AG93-Lem2} and obtained a generalization of the G\"{o}llnitz theorem,
a higher level extension of Schur's partition theorem. Alladi, Andrews and
Berkovich \cite[Eq.\,\,(1.7)]{AAB} further obtained a remarkable four parameter extension of
the identity \eqref{AG93-Lem2}, which implies a four parameter
partition theorem and thereby extends the G\"{o}llnitz theorem.

Due to its significance, the Alladi-Gordon key identity certainly deserves to be further studied. Alladi and Gordon \cite{Alladi-Gordon-93} gave three different proofs of \eqref{AG93-Lem2} by using weighted words, $q$-Chu-Vandermonde summation formula and the inclusion-exclusion principal respectively. In this paper we will present an overpartition interpretation of the left-hand side of \eqref{eq-ij2} and then give a combinatorial proof of the Alladi-Gordon key identity.

This paper is organized as follows. In Section 2 we will review Joichi-Stanton's insertion algorithm for partitions and then give an overpartition interpretation of the Alladi-Gordon key identity. In Section 3, based on this interpretation, we will establish an involution and give a combinatorial proof of the Alladi-Gordon key identity.

\section{An overpartition interpretation of the Alladi-Gordon key identity}\label{sec2}

The aim of this section is to give a combinatorial interpretation of the left-hand side of
\eqref{eq-ij2} in terms of overpartitions. This is achieved by applying Joichi-Stanton's insertion algorithm for partitions.

Let us first review some definitions and notations about partitions.
Recall that a \emph{partition} $\lambda$ of $n\in\mathbb{N}$ with $k$ parts is
denoted by a vector $\lambda=(\lambda_1, \lambda_2, \ldots,
\lambda_k)$, where $\lambda _1 \geq  \lambda_2 \geq \cdots \geq
\lambda_k \geq 0$ and $\sum_{i=1}^k \lambda_i=n$. The number $n$ is called the \emph{size} of $\lambda$, denoted by
$|\lambda|$. For convenience, the \emph{length} of $\lambda$ is defined to be the number $k$ of nonnegative parts of
$\lambda$, denoted by $\ell(\lambda)$. (Note that $\ell(\lambda)$ usually enumerates the number of positive parts.)
An \emph{overpartition}
is a partition in which the first occurrence of a number may be
overlined. For example, $\lambda=(9,\overline{7},6,5,
5,\overline{2},2,\overline{1})$ is an overpartition with three
overlined parts. An ordinary partition can also be treated as an
overpartition with no overlined parts. The concept of overpartition was first proposed by Corteel and
Lovejoy\cite{CorLov2004} while studying basic hypergeometric series.
For deeper researches on overpartitons, see \cite{LovMal2008, BriLov2009, Sills2010}.

An overpartition can also be
understood as a pair of partitions $(\alpha, \beta)$, where $\alpha$
is a partition with distinct parts and $\beta$ is an ordinary
partition. Joichi and Stanton \cite{JoiSta} found the following
fundamental bijection which can be restated in terms of overpartitions.

\begin{theo}\label{theo2.1} There is a one-to-one correspondence between
overpartitions with $n$ nonnegative parts, and pairs of partitions
$(\alpha, \beta)$, where $\alpha$ is a partition with distinct parts
from the set $\{0, 1, 2, \ldots, n-1\}$ and $\beta$ is a partition with $n$
nonnegative parts.
\end{theo}

The above correspondence can be described as an insertion algorithm
\cite[Algorithm $\Phi$]{JoiSta}. Given an ordinary partition
$\beta$, we may insert a part $m$ into $\beta$, by adding 1 to the
first $m$ parts of $\beta$, and putting an overline above the
$(m+1)$-th part. Moreover, we can add other distinct parts in the
same way.
\begin{exam}If $\alpha=(5,3,0)$ and $\beta=(9,6,5,2,2,0)$, then we
get an overpartition
$(\overline{11},8,7,\overline{3},3,\overline{0})$.
\end{exam}

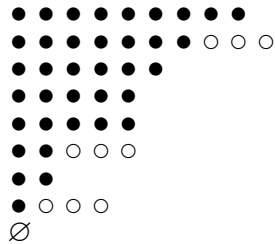
\begin{figure}[h]
\vskip -5pt
\setlength\unitlength{5.2pt}
\begin{center}
\begin{picture}(26,20)
\multiput(6.0,17.0)(2,0){9}{\circle*{1}}
\multiput(6.0,15.0)(2,0){7}{\circle*{1}}
\multiput(20,15)(2,0){3}{\circle{1}}
\multiput(6.0,13.0)(2,0){6}{\circle*{1}}
\multiput(6.0,11.0)(2,0){5}{\circle*{1}}
\multiput(6.0,9.0)(2,0){5}{\circle*{1}}
\multiput(6.0,7.0)(2,0){2}{\circle*{1}}
\multiput(10.0,7.0)(2,0){3}{\circle{1}}
\multiput(6.0,5.0)(2,0){2}{\circle*{1}}
\put(6.0,3.0){\circle*{1}}
\multiput(8.0,3.0)(2,0){3}{\circle{1}}
\put(5.2,0.5){$\varnothing$}
\end{picture}
\caption{An overpartition $\lambda=(9,\overline{7},6,5,
5,\overline{2},2,\overline{1},0)$ with a weight $3$ in each overline. \label{f-3-1}}
\end{center}
\vskip -5mm
\end{figure}

To give a combinatorial interpretation of \eqref{eq-ij2}, we shall assign a weight to
each overlined part of an overpartition. As in \cite{ChenZhao2005}, each overlined part of an overpartition has the same weight. For example, $\lambda=(9,\overline{7},6,5,
5,\overline{2},2,\overline{1},0)$ with a weight $3$ endowed in each overline
is displayed in Figure \ref{f-3-1}, where each overline of $\lambda$
is represented by a row of three hollow dots, and the part zero is represented by the symbol $\varnothing$.

Given $0\leq k\leq j\leq i$, let $A(i,k)$ denote the set of partitions into distinct parts
from the set $\{i-k+1,i-k+2,\ldots,i\}$ plus the empty partition, and let $B(j,k)$ denote the set of partitions into
$k$ nonnegative parts with each part not exceeding $j-k$. It is well known that
\begin{align}\label{eq-gauss}
\sum_{\lambda\in A(i,k)}(-1)^{\ell(\lambda)}q^{|\lambda|}=(q^{i-k+1};q)_k,\quad\quad
\sum_{\lambda\in B(j,k)}q^{|\lambda|}={j\brack k}.
\end{align}

We now come to the main result of this section. For a weighted overpartition $\lambda$, let $ol(\lambda)$ denote the number of overlined parts of $\lambda$, and let $w(\lambda)$ denote the weight assigned to each overlined part of $\lambda$. Given $0\leq k\leq j\leq i$,
let $O(i,j,k)$ denote the set of weighted overpartitions into $k$ nonnegative parts, which satisfy the following three conditions:
\begin{itemize}
\item[(1)] each part is less than or equal to $j-1$;

\item[(2)] for $k\geq 2$ and each $1\leq s\leq k-1$ there are at least $k-s$ overlined parts to the right of $j-s$ if it occurs as a part;

\item[(3)] and each overline is endowed with a weight $i-k+1$.
\end{itemize}
Note that the empty partition is the sole element of $O(i,j,0)$.
For fixed $i,j$ satisfying $0\leq j\leq i$, let
\begin{align}\label{eq-set}
O(i,j)=\biguplus\limits_{k=0}^j O(i,j,k).
\end{align}
The main result of this section is as follows.

\begin{theo}\label{theo2.2} Given $0\leq j\leq i$, we have
\begin{align}
\sum\limits_{k=0}^{j}(q^{i-k+1};q)_k\,{j\brack k} q^{(i-k)(j-k)}=\sum_{\lambda\in O(i,j)}(-1)^{ol(\lambda)}q^{|\lambda|+ol(\lambda)w(\lambda)}q^{(i-\ell(\lambda))(j-\ell(\lambda))}.
\end{align}
\end{theo}

\pf By \eqref{eq-set}, it suffices to prove that
$$(q^{i-k+1};q)_k\,{j\brack k}=\sum_{\lambda\in O(i,j,k)}(-1)^{ol(\lambda)}q^{|\lambda|+ol(\lambda)(i-k+1)},$$
since for each $\lambda\in O(i,j,k)$ we have $\ell(\lambda)=k$ and $w(\lambda)=i-k+1$.
It is clear true for $k=0$. In the following we may assume that $k\geq 1$.
In view of \eqref{eq-gauss}, we only need to give a weight-preserving bijection $\overline{\Phi}$
between the set $A(i,k)\times B(j,k)$ and the set $O(i,j,k)$. Actually, we can obtain $\overline{\Phi}$
by using Joichi-Stanton's insertion algorithm.

For any given pair $(\gamma,\beta)\in A(i,k)\times
B(j,k)$, define $\overline{\Phi}(\gamma,\beta)$ to be the partition $\lambda$ obtained as follows.

\begin{itemize}
\item[(i)] If $\gamma$ is the empty partition, then let $\lambda=\beta$. By Property (2) of the definition of $O(i,j,k)$, it is clear that $B(j,k)\subseteq O(i,j,k)$. Therefore, in this case we have $\lambda\in O(i,j,k)$.

\item[(ii)] If $\gamma$ is not the empty partition, then let $\overline{\gamma}$ denote the partition obtained from $\gamma$ by decreasing each part by $i-k+1$. Therefore, $\overline{\gamma}$ is a partition into distinct parts from the set $\{0,1,\ldots,k-1\}$. Now we insert $\overline{\gamma}$
into $\beta$ by applying Joichi-Stanton's insertion algorithm, and obtain an overpartition $\lambda$ with at least one overlined parts. If each overline is endowed with a weight $i-k+1$, then it is routine to verify that the weighted overpartition $\lambda$ lies in $O(i,j,k)$. Note that the number of parts of $\gamma$ is equal to the number of overlined parts of $\lambda$. Thus
$$(-1)^{\ell(\gamma)}q^{|\gamma|}q^{|\beta|}=(-1)^{ol(\lambda)}q^{|\lambda|+ol(\lambda)(i-k+1)}.$$
\end{itemize}

It remains to show that $\overline{\Phi}$ is reversible. There are two cases to consider.
\begin{itemize}

\item[(i')] If $\lambda\in O(i,j,k)$ and there are no overlined parts in $\lambda$, then again by Property (2) of the definition of $O(i,j,k)$, we must have $\lambda\in B(j,k)$. In this case, let $\overline{\Phi}^{-1}(\lambda)=(\gamma,\beta)$, where $\gamma$ is the empty partition and $\beta=\lambda$.

\item[(ii')] If $\lambda\in O(i,j,k)$ and there are at least one overlined parts in $\lambda$, then by reversing the insertion algorithm, we will obtain a pair of partitions $(\overline{\gamma},\beta)$. Clearly, $\overline{\gamma}$ is a partition into distinct parts from the set $\{0,1,\ldots,k-1\}$ since there are $k$ parts in $\lambda$. It is also clear that $\beta$ has only $k$ parts. We further need to show that each part of $\beta$ is not exceeding $j-k$. Suppose that there are $t$ overlined parts to the right of $\lambda_1$, then $\beta_1=\lambda_1-t$. Assume that $\lambda_1=j-s$ for some $1\leq s\leq k-1$. By Property (2) of the definition of $O(i,j,k)$ we have $t\geq k-s$. Therefore, $\beta_1=\lambda_1-t=j-s-t\leq j-k$.

\end{itemize}

This completes the proof.
\qed

The following example gives an illustration of the map $\overline{\Phi}$ of the above proof.

\begin{exam}\label{exam3.1}
For $i=9,\ j=6,\ k=3$, $\gamma=(8,7)\in A(9,3)$, $\beta=(3,3,2)\in
B(6,3)$. We shall transform $(\gamma,\beta)$ into
$\lambda\in O(9,6,3)$ in two steps.
\begin{itemize}
\item[(1)] Change $\gamma=(8,7)$ into
$\overline{\gamma}=({1},{0})$ by decreasing each part by $7$.

\item[(2)] Insert $\overline{\gamma}=({1},{0})$ into $\beta=(3,3,2)$ to
obtain an overpartition $\lambda=(\overline{4},\overline{3},2)\in
O(9,6,3)$, where each overline contains a weight $7$. See
Figure \ref{f-Fg-1}.
\end{itemize}
By reversing the procedure it is easy to obtain $(\gamma,\beta)$
from $\lambda$.
\end{exam}

\begin{center}
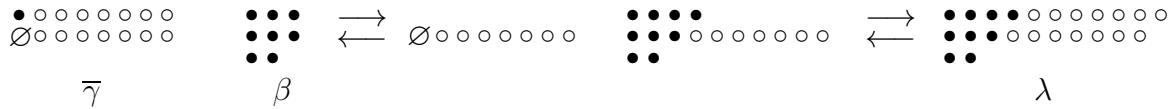
\begin{figure}[hbt]
\setlength\unitlength{4.0pt}
\begin{picture}(110,10)
\put(2.0,9.0){\circle*{1}}
\multiput(4.0,9.0)(2,0){7}{\circle{1}}
\put(1.0,6.1){$\varnothing$}
\multiput(4.0,7.0)(2,0){7}{\circle{1}}
\put(8.0,1.0){$\overline{\gamma}$}
\multiput(24.0,9.0)(2,0){3}{\circle*{1}}
\multiput(24.0,7.0)(2,0){3}{\circle*{1}}
\multiput(24.0,5.0)(2,0){2}{\circle*{1}}
\put(26.0,1.0){$\beta$}
\put(32.0,8.0){$\longrightarrow$}
\put(32.0,6.0){$\longleftarrow$}
\put(38.7,6.2){$\varnothing$}
\multiput(42.0,7.0)(2,0){7}{\circle{1}}
\multiput(60.0,9.0)(2,0){4}{\circle*{1}}
\multiput(60.0,7.0)(2,0){3}{\circle*{1}}
\multiput(66.0,7.0)(2,0){7}{\circle{1}}
\multiput(60.0,5.0)(2,0){2}{\circle*{1}}
\put(82.0,8.0){$\longrightarrow$}
\put(82.0,6.0){$\longleftarrow$}
\multiput(90.0,9.0)(2,0){4}{\circle*{1}}
\multiput(98.0,9.0)(2,0){7}{\circle{1}}
\multiput(90.0,7.0)(2,0){3}{\circle*{1}}
\multiput(96.0,7.0)(2,0){7}{\circle{1}}
\multiput(90.0,5.0)(2,0){2}{\circle*{1}}
\put(98.0,1.0){$\lambda$}
\end{picture}
\caption{Insertion of
$\overline{\gamma}=({1},{0})$ into $\beta=(3,3,2)$
leads to $\lambda=(\overline{4},\overline{3},2)$, where $\overline{\gamma}$ is represented as the  overpartition $(\overline{1},\overline{0})$ and each
overline has a weight $7$ endowed. \label{f-Fg-1}}
\end{figure}
\end{center}

\section{Combinatorial proof of the Alladi-Gordon key identity}\label{secSchurInv}

The aim of this section is to prove the following result by constructing an involution on the set $O(i,j)$.
\begin{theo}\label{theo3.1} Given $0\leq j\leq i$, we have
\begin{align}\label{eq-3.1}
\sum_{\lambda\in O(i,j)}(-1)^{ol(\lambda)}q^{|\lambda|+ol(\lambda)w(\lambda)}q^{(i-\ell(\lambda))(j-\ell(\lambda))}=1,
\end{align}
where $O(i,j)$ is as defined in \eqref{eq-set}.
\end{theo}
Combining Theorem \ref{theo2.2}, this provides a combinatorial proof of the the Alladi-Gordon key identity.

To prove Theorem \ref{theo3.1}, we first give a decomposition of $O(i,j)$. For $\lambda\in O(i,j)$ let $\overline{\lambda_{t}}$ denote the largest overlined part of
$\lambda$. Let
\begin{align*}
{O}_1(i,j)&=\left\{\left.\lambda \in O(i,j)\, \right|\,
ol(\lambda)\geq1, \lambda_t=j-\ell(\lambda)+ol(\lambda)-1 \right\},\\[5pt]
{O}_2(i,j)&=\left\{\left.\lambda \in O(i,j)\, \right|\,
ol(\lambda)\geq 1, \lambda_t<j-\ell(\lambda)+ol(\lambda)-1 \right\},\\[5pt]
{O}_3(i,j)&=\left\{\left.\lambda \in O(i,j)\, \right|\,
ol(\lambda)=0 \right\}.
\end{align*}
For the convenience, let the empty partition belong to ${O}_3(i,j)$.

\begin{lem} For $0\leq j\leq i$, we have
$$O(i,j)={O}_1(i,j) \uplus {O}_2(i,j) \uplus {O}_3(i,j).$$
\end{lem}
\pf It is clear that ${O}_1(i,j),{O}_2(i,j)$ and ${O}_3(i,j)$ are disjoint from each other. It suffices to show that for each $\lambda\in O(i,j)$ with $ol(\lambda)\geq 1$, we have $\lambda_t\leq j-\ell(\lambda)+ol(\lambda)-1$.
Otherwise, suppose that $\lambda_t=j-s$ and $s<\ell(\lambda)-ol(\lambda)+1$. By Property (2) of the definition of $O(i,j,k)$, there are at least $\ell(\lambda)-s\geq ol(\lambda)$ overlined parts to the right of $\lambda_t$, contradicting the definition of $ol(\lambda)$. This completes the proof.
\qed

With the above decomposition of $O(i,j)$, now we can give a bijective proof of Theorem \ref{theo3.1}.

\noindent \textit{Proof of Theorem \ref{theo3.1}.}
For $\lambda\in O(i,j)$, let
$$f(\lambda)=(-1)^{ol(\lambda)}q^{|\lambda|+ol(\lambda)w(\lambda)}q^{(i-\ell(\lambda))(j-\ell(\lambda))}.$$
To give a bijective proof, we define an involution, denoted $\Psi$, acting on $O(i,j)$ as follows:

\begin{itemize}
\item[{(1)}] If $\lambda\in O_1(i,j)$, then let $\Psi(\lambda)$ denote the overpartition obtained from $\lambda$ by removing the
largest overlined part $\lambda_t$. In this case, we have
\begin{align*}
\ell(\Psi(\lambda))&=\ell(\lambda)-1\\
ol(\Psi(\lambda))&=ol(\lambda)-1\\
|\Psi(\lambda)|&=|\lambda|-\lambda_t\\
&=|\lambda|-(j-\ell(\lambda)+ol(\lambda)-1).
\end{align*}
By Property 3 of the definition of $O(i,j,k)$, we have
$$w(\Psi(\lambda))=i-\ell(\Psi(\lambda))+1=(i-\ell(\lambda)+1)+1=w(\lambda)+1$$
It is routine to verify that $f(\Psi(\lambda))+f(\lambda)=0$. Note that
if $ol(\lambda)=1$, then clearly $\Psi(\lambda)\in O_3(i,j)$. If $ol(\lambda)>1$,
then we must have $\Psi(\lambda)\in O_2(i,j)$ since
$$\Psi(\lambda)_t<\lambda_t=j-\ell(\lambda)+ol(\lambda)-1=j-\ell(\Psi(\lambda))+ol(\Psi(\lambda))-1.$$
(The inequality $\Psi(\lambda)_t<\lambda_t$ follows from the definition of overpartitions.)
In both cases, we have $\ell(\Psi(\lambda))=\ell(\lambda)-1\leq j-1$.

\item[{(2)}] If $\lambda\in O_2(i,j)$, then we must have $j>\ell(\lambda)$. Otherwise, if $j=\ell(\lambda)$, then $\lambda_t<ol(\lambda)-1$, contradicting the fact that the overlined parts of an overpartition are distinct nonnegative integers. Then let $\Psi(\lambda)$ denote the overpartition obtained from $\lambda$ by inserting an overlined part $j-\ell(\lambda)+ol(\lambda)-1$. Clearly,
    $\Psi(\lambda)\in O_1(i,j)$ and $\Psi(\Psi(\lambda))=\lambda$.

\item[{(3)}] If $\lambda\in O_3(i,j)$, then we define $\Psi(\lambda)$ as follows according to whether $j>\ell(\lambda)$. If $j>\ell(\lambda)$, then let $\Psi(\lambda)$ denote the overpartition obtained from $\lambda$ by inserting an overlined part $j-\ell(\lambda)-1$. In this case, it is clear that  $\Psi(\lambda)\in O_1(i,j)$ and $\Psi(\Psi(\lambda))=\lambda$. If $j=\ell(\lambda)$, then $\lambda$ must be the partition $(\underbrace{0,0,\ldots,0}_{j's})$. Otherwise,  we will have $\lambda_1>0$, and by Property 2 of $O(i,j,j)$ there must be at least one overlined part in $\lambda$ contradicting $ol(\lambda)=0$.
    In this case let $\Psi(\lambda)=\lambda$.
\end{itemize}
By the involution $\Psi$ of $O(i,j)$, we have
$$\sum_{\lambda\in O(i,j)} f(\lambda)=\sum_{\lambda=(\underbrace{0,0,\ldots,0}_{j's})} f(\lambda)=1.$$
This completes the proof. \qed

In fact, there is a graphical representation of the involution $\Psi$ of $O(i,j)$ in the above proof, which seems more convenient and intuitive. Since each $\lambda\in O(i,j)$ contributes a term
$$f(\lambda)=(-1)^{ol(\lambda)}q^{|\lambda|+ol(\lambda)w(\lambda)}q^{(i-\ell(\lambda))(j-\ell(\lambda))},$$
we may consider $\lambda$ as a pair of partitions $(\lambda, \hat{\lambda})$, where $\hat{\lambda}$ is the unique rectangular partition $(\underbrace{i-\ell(\lambda),\ldots,i-\ell(\lambda)}_{(j-\ell(\lambda))'s})$.

\begin{exam} Take $i=9,\, j=6$ and let $\lambda=(\overline{4},\overline{3},2)$. In this case, we have $\hat{\lambda}=(6,6,6)$,
$ol(\lambda)=2$, $\ell(\lambda)=3$, $w(\lambda)=i-\ell(\lambda)+1=7$, and hence $\lambda\in O_1(i,j)$. Thus, $\Psi(\lambda)=(\overline{3},2)\in O_2(i,j)$ and $\hat{\Psi(\lambda)}=(7,7,7,7)$. Geometrically, $\Psi$ acts on $\lambda$ (or equivalently $(\lambda,\hat{\lambda})$) as illustrated in Figure \ref{f-Fg-2}: remove a row of dots representing the largest overlined part, add a hollow dot to the rightmost of each overlined part, and
append a hook to the top-left of the diagram of $\hat{\lambda}$.

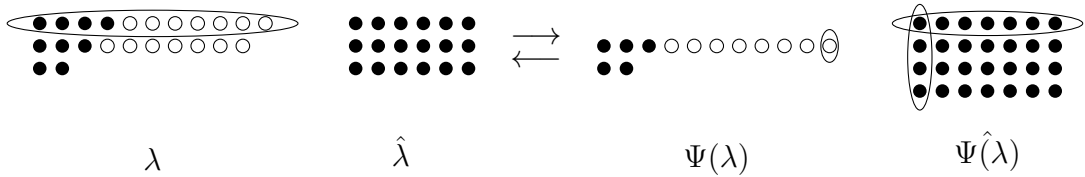
\begin{figure}[ht]
\begin{center}
\begin{tikzpicture}
\foreach \x in {-1.5,-1.2,...,-0.31}
\fill (\x,0) circle (0.9mm);
\foreach \x in {-0.3,0,...,1.51}
\draw (\x,0) circle (0.9mm);
\foreach \x in {-1.5,-1.2,-0.9}
\fill (\x,-0.3) circle (0.9mm);
\foreach \x in {-0.6,-0.3,...,1.21}
\draw (\x,-0.3) circle (0.9mm);
\fill (-1.5,-0.6)  circle(0.9mm);
\fill (-1.2,-0.6)  circle(0.9mm);
\draw (0,0)  ellipse (55pt and 5.0pt);
\draw(0,-1.8) node {$\lambda$};
\foreach \x in {2.7,3.0,...,4.21}
\fill (\x,0) circle (0.9mm);
\foreach \x in {2.7,3.0,...,4.21}
\fill (\x,-0.3) circle (0.9mm);
\foreach \x in {2.7,3.0,...,4.21}
\fill (\x,-0.6) circle (0.9mm);
\draw(3.3,-1.7) node {$\hat{\lambda}$};
\draw(5.1,-0.2) node {$\longrightarrow$};
\draw(5.1,-0.5) node {$\longleftarrow$};
\foreach \x in {6.0,6.3,6.6}
\fill (\x,-0.3) circle (0.9mm);
\foreach \x in {6.9,7.2,...,9.1}
\draw (\x,-0.3) circle (0.9mm);
\foreach \x in {6.0,6.3}
\fill (\x,-0.6) circle (0.9mm);
\draw (0,0) (9.0,-0.3)  [rotate=0] ellipse (3.3pt and 6.2pt);
\draw(7.5,-1.8) node {$\Psi(\lambda)$};
\foreach \x in {10.2,10.5,...,12.1}
\fill (\x,0) circle (0.9mm);
\foreach \x in {10.2,10.5,...,12.1}
\fill (\x,-0.3) circle (0.9mm);
\foreach \x in {10.2,10.5,...,12.1}
\fill (\x,-0.6) circle (0.9mm);
\foreach \x in {10.2,10.5,...,12.1}
\fill (\x,-0.9) circle (0.9mm);
\draw (0,0) (11.1,0) ellipse (36pt and 4.6pt);
\draw (0,0) (10.2,-0.45) ellipse (4.6pt and 20pt);
\draw(11.1,-1.7) node {$\hat{\Psi(\lambda)}$};
\end{tikzpicture}
\end{center}
\caption{The involution $\Psi$.} \label{f-Fg-2}
\end{figure}

\end{exam}

\noindent{\bf Acknowledgments.}
I would like to thank Arthur L.B. Yang, Qing-Hu Hou, and Guoce Xin
for valuable comments leading to an improvement of an earlier
version. This work was supported by the National Natural Science
Foundation of China Grants No.\,70871011 and No.\,71171018, Program
for New Century Excellent Talents in University, and the Fundamental
Research Funds for the Central Universities.

\end{document}